\documentclass[12pt]{article}
\usepackage{amsfonts,amsmath,amssymb,amscd}
\newtheorem{definition}{Definition}[section]

\newtheorem{proposition}[definition]{Proposition}
\newtheorem{theorem}[definition]{Theorem}
\newtheorem{corollary}[definition]{Corollary}
\newtheorem{lemma}[definition]{Lemma}

\def\cA{{\cal A}}                    
\def\cD{{\cal D}}          \def\cE{{\cal E}}          \def\cF{{\cal F}}
\def\cG{{\cal G}}          \def\cH{{\cal H}}

\def\fS{{\mathfrak S}}

\def\fa{{\mathfrak a}}
\def\fg{{\mathfrak g}}
\def\fl{{\mathfrak l}}
\def\fs{{\mathfrak s}}

\newcommand{\CC}{{\mathbb C}}

\newcommand{\ZZ}{{\mathbb Z}}
\newcommand{\eps}{{\varepsilon}}

\newcommand{\finproof}{{\hfill \rule{5pt}{5pt}}}

\def\QTHA{Quasitriangular Hopf Algebra (QTHA)\def\QTHA{QTHA}}
\def\QTQHA{Quasitriangular Quasi-Hopf Algebra (QTQHA)\def\QTQHA{QTQHA}}

\begin{document}
\pagestyle{empty}


\begin{center}

{\Large \textsf{Hopf structure of the Yangian
$Y(sl_n)$ in the Drinfel'd realisation}} 

\vspace{10mm}

{\large N.~Cramp{\'e}}

\vspace{10mm}

\emph{Laboratoire d'Annecy-le-Vieux de Physique Th{\'e}orique}

\emph{LAPTH, CNRS, UMR 5108, Universit{\'e} de Savoie}

\emph{B.P. 110, F-74941 Annecy-le-Vieux Cedex, France}

\vspace{7mm}

\end{center}

\vfill
\vfill

\begin{abstract}
The Yangian of the Lie algebra $sl_n$ is known to have different
presentations, in particular the RTT realisation and
the Drinfel'd realisation.
Using the isomorphism between them, 
the explicit expressions of the 
comultiplication, the antipode and the counit in the Drinfel'd realisation 
of the Yangian $Y(sl_n)$ are given. As examples, the cases of $Y(sl_2)$
and $Y(sl_3)$ are worked out.
\end{abstract}

\vfill
\vfill

\rightline{LAPTH-978/03}
\rightline{math.QA/0304254}

\baselineskip=16pt

\newpage

\pagestyle{plain}
\setcounter{page}{1}
The Yangian $Y(\fa)$ based on a simple Lie algebra $\fa$ is defined
\cite{Dri85,Dri86} as the
(unique) homogeneous quantisation of the algebra $\fa[u] = \fa \otimes
\CC[u]$ endowed with its standard bialgebra structure, where $\CC[u]$
is the ring of polynomials in the indeterminate $u$. This algebra has
a structure of a non-cocomutative Hopf algebra, which partially 
explains the importance of Yangians and their representations 
in the study of quantum inverse problem. Among the different
presentations of the Yangians, the one known as the 
Drinfel'd realisation is well adapted for the study 
of their representations \cite{Dri88}. No explicit
formula for the Hopf structure in this realisation was known yet,
except for $sl_2$ \cite{Molev2001} and for $osp(1|2)$ \cite{yop}. The
aim of this letter is to give an explicit expression of the
comultiplication, the antipode and the counit in the Drinfel'd 
realisation for $Y(sl_n)$.
Note that partial results were given in \cite{CP,KT}.
The comultiplication given in this letter can be extended to the
double Yangian $DY(sl_n)$. One can show that the so-called Drinfel'd
comultiplication defined only for the double
Yangian is the twist of this extended comultiplication
(see for example \cite{KT,Iohara}).\\

This letter is organised as follows. In the first section,
the RTT formalism \cite{FRT} and Drinfel'd realisation of $Y(sl_n)$  
are presented, which allow us to give the normalisation of
the generators as well as the exact form of the $R$-matrix and 
of the quantum determinant.
In the second section, some properties about the quantum minors,
needed in the following,
are explained. The expressions of 
the isomorphism using the quantum minors or the Gauss decomposition
are then
presented. The main theorem of this letter, i.e.
the explicit form of the Hopf structure, is exposed in the next two
sections. Finally, as illustrative examples, the $Y(sl_2)$ and
$Y(sl_3)$ cases are worked out.

\section{Two realisations of the Yangian $Y(sl_n)$}
\setcounter{equation}{0}

In this section, two different realisations of the Yangian based on
the Lie algebra $sl_n$ are presented: the RTT formalism and the
Drinfel'd realisation \cite{Dri88}. \\
The first realisation uses the RTT formalism
 \cite{Dri88,Molev2001,FRT}.
Let $V^{(n)}$ denotes the $n$-dimensional fundamental vector space representation
of $sl_n$. 
The Yang's R-matrix is given by 
\begin{eqnarray}
&&R^{(n)}_{12}(u)=I\otimes I+\sum_{1\le i,j \le n}\frac{E_{ij} \otimes
E_{ji}}{u} \in End(V^{(n)}\otimes V^{(n)})
\end{eqnarray} 
where $E_{ij}$ is the elementary matrix with entry 1 in row $i$ and column $j$ and 0
elsewhere.
This R-matrix satisfies the following properties
\begin{eqnarray}
&&R^{(n)}_{12}(u)R^{(n)}_{13}(u+v)R^{(n)}_{23}(v) 
=R^{(n)}_{23}(v)R^{(n)}_{13}(u+v)R^{(n)}_{12}(u)\quad\mbox{(Yang-Baxter
equation)}\\
&&R^{(n)}_{12}(u)R^{(n)}_{21}(-u)=\frac{u^2-1}{u^2}\;(I\otimes I)\quad\mbox{(unitarity)}
\end{eqnarray}

\begin{theorem}
The Yangian of $sl_n$, $Y(sl_n)$, is isomorphic to the associative algebra,
$U(R)$, generated by the unit and the elements $\{T_{i,j}^{(k)}\;|\;1\le i,j
\le n,k\in\ZZ_{>0}\}$ gathered in the formal series
\begin{eqnarray} 
T(u) 
  =1 + \sum_{i,j=1}^{n} \sum_{n \in \ZZ_{> 0}}
  T_{i,j}^{(n)} \, u^{-n} \, E_{ij}
  = \sum_{i,j=1}^{n} {T}_{i,j}(u) \, E_{ij}
\end{eqnarray}
subject to the defining relations
\begin{eqnarray}
\label{RTT}
&&R^{(n)}_{12}(u-v)(T(u)\otimes 1)(1\otimes T(v))
=(1\otimes T(v))(T(u)\otimes 1)R^{(n)}_{12}(u-v)\quad\\
\mbox{and}&&qdet\,T(u)=1
\end{eqnarray}
where
$\quad qdet\,T(u)=\sum_{\sigma\in\mathfrak{S}_n}sgn(\sigma)~
T_{\sigma(1),1}(u)\cdots T_{\sigma(n),n}(u+n-1)$.
\end{theorem}
The defining relations (\ref{RTT}) as commutators of
${T}_{ij}(u)$:
\begin{eqnarray}
\label{RTT_relcom}
&&-(u-v)\left[{T}_{i,j}(u),{T}_{k,l}(v)\right]
={T}_{k,j}(u)\,{T}_{i,l}(v)-{T}_{k,j}(v)\,{T}_{i,l}(u).
\end{eqnarray}
The map 
\begin{eqnarray}
\label{inverse}
T(u)&\longmapsto& T^{-1}(-u)\equiv T^*(u)
\end{eqnarray}
defines an automorphism of $Y(sl_n)$.\\
To avoid ambiguity, let us stress that
$qdet\,T(-u)=\sum_{\sigma\in\mathfrak{S}_n}sgn(\sigma)~
T_{\sigma(1),1}(-u)\cdots T_{\sigma(n),n}(-u+n-1)$ is different from
the quantum determinant of the matrix $\widetilde{T}(u)=T(-u)$:
\begin{equation} 
qdet\,\widetilde{T}(u)=\sum_{\sigma\in\mathfrak{S}_n}sgn(\sigma)~
T_{\sigma(1),1}(-u)\cdots T_{\sigma(n),n}(-u-n+1).
\end{equation} 
The Yangian $Y(sl_n)$ has a Hopf algebra structure and 
the explicit forms of comultiplication, 
antipode and counit are
\begin{eqnarray}
\label{hopf}
&&\Delta(T_{i,j}(u))=\sum_{k=1}^{n}T_{i,k}\otimes T_{k,j}\;,\qquad
S(T(u))=T^{-1}(u)
\qquad\mbox{and}\qquad\epsilon(T_{i,j}(u))=\delta_{ij}\;.
\end{eqnarray}
The second realisation of the Yangian uses the so-called Drinfel'd generators.
Let $\{\alpha_i|1\le i \le n-1\}$ be the set of simple
roots of $sl_n$ and $(\cdot,\cdot)$ be the standard non-degenerate symmetric
invariant bilinear form on $sl_n$. For each simple root
$\alpha_i$, $e_{\alpha_i}$ and $f_{\alpha_i}$ are the corresponding root vectors, such that
$(e_{\alpha_i},f_{\alpha_i})=1$, and
$h_{\alpha_i}=[e_{\alpha_i},f_{\alpha_i}]$ are the Cartan generators.
The Drinfel'd realisation of the Yangian is given by the following
theorem \cite{Dri88}.
\begin{theorem}
The Yangian of $sl_n$, $Y(sl_n)$, is isomorphic to the associative algebra
$\cA$, generated by the unit and the elements 
$\{e^{(k)}_{i},f^{(k)}_{i},h_{i}^{(k)}\,|\,1\le i \le
n-1,\, k\in \ZZ_{\ge 0}\}$ subject to the defining relations 
\begin{eqnarray}
&&[h_{i}^{(k)},h_{j}^{(l)}]=0\;,\qquad
[e^{(k)}_{i},f^{(l)}_{j}]=\delta_{i,j}~h^{(k+l)}_{i}\;,\\
&&[h_{i}^{(0)},e^{(l)}_{j}]=(\alpha_i,\alpha_j)~e^{(l)}_{j}\;,\qquad
[h_{i}^{(0)},f^{(l)}_{j}]=-(\alpha_i,\alpha_j)~f^{(l)}_{j}\;,\\
&&[h_{i}^{(k+1)},e^{(l)}_{j}]-[h_{i}^{(k)},e^{(l+1)}_{j}] 
=\frac{1}{2}~(\alpha_i,\alpha_j)~
(h_{i}^{(k)}\,e^{(l)}_{j}+e^{(l)}_{j}\,h^{(k)}_{i})\;,\\
&&[h_{i}^{(k+1)},f^{(l)}_{j}]-[h_{i}^{(k)},f^{(l+1)}_{j}] 
=-\frac{1}{2}~(\alpha_i,\alpha_j)~
(h_{i}^{(k)}\,f^{(l)}_{j}+f^{(l)}_{j}\,h^{(k)}_{i})\;,\\
&&[e^{(k+1)}_{i},e^{(l)}_{j}]-[e^{(k)}_{i},e^{(l+1)}_{j}] 
=\frac{1}{2}~(\alpha_i,\alpha_j)~
(e^{(k)}_{i}\,e^{(l)}_{j}+e^{(l)}_{j}\,e^{(k)}_{i})\;,\\
&&[f^{(k+1)}_{i},f^{(l)}_{j}]-[f^{(k)}_{i},f^{(l+1)}_{j}] 
=-\frac{1}{2}~(\alpha_i,\alpha_j)~
(f^{(k)}_{i}\,f^{(l)}_{j}+f^{(l)}_{j}\,f^{(k)}_{i})\;,
\end{eqnarray} 
and to the Serre relations, for 
$i\neq j$ and $n_{ij}=1-2\frac{(\alpha_i,\alpha_j)}
{(\alpha_i,\alpha_i)}$:
\begin{eqnarray} 
&&\sum_{\sigma\in\mathfrak{S}_{n_{ij}}}[e^{(k_{\sigma(1)})}_{i},[\cdots,
[e^{(k_{\sigma(n_{ij})})}_{i},e^{(l)}_{j}]\cdots]=0\;,\\
&&
\sum_{\sigma\in\mathfrak{S}_{n_{ij}}}[f^{(k_{\sigma(1)})}_{i},[\cdots,
[f^{(k_{\sigma(n_{ij})})}_{i},f^{(l)}_{j}]\cdots]=0\;.
\end{eqnarray} 
\end{theorem}
For later conveniences, we define the following formal series:
\begin{eqnarray}
e_{i}(u)=\sum_{k=0}^{+\infty}\frac{{e^{(k)}_{i}}}{u^{k+1}}
,~
f_{i}(u)=\sum_{k=0}^{+\infty}\frac{{f^{(k)}_{i}}}{u^{k+1}}
~~\mbox{and}~~
h_i(u)=1+\sum_{k=0}^{+\infty}\frac{h_{i}^{(k)}}{u^{k+1}}
&,&\mbox{for}\quad 1\le i \le n-1\;. 
\end{eqnarray} 
The mapping $e_{\alpha_i}\longmapsto e^{(0)}_{i}\;,
\;f_{\alpha_i}\longmapsto f^{(0)}_{i}\;,
\;h_{\alpha_i}\longmapsto h^{(0)}_{i}$ defines an
embedding $U(sl_n)\hookrightarrow Y(sl_n)$, where $U(sl_n)$ is the
universal enveloping algebra of $sl_n$. 
\\
\section{Quantum minors}
\setcounter{equation}{0}
Before giving the expression of the
isomorphism that relates the two Yangian presentations, in the next
section, we introduce the notion of quantum minors and give 
some of their properties.\\
Let $I=\{a_1,a_2,\ldots,a_m\}$ and
$J=\{b_1,b_2,\ldots,b_m\}$ such that $I,J\subset\{1,\ldots,n\}$ and  
$card(I)=card(J)=m$ with  $1\le m\le n$.
The set of generators $\{T_{a_i,b_j}(u)|1\le i,j \le m\}$ defines a
subalgebra of $Y(sl_n)$ with the following commutation relations
\begin{eqnarray}
&&\hspace{-20mm}R^{(m)}_{12}(u-v)
\left(T\!\left(^{a_1\cdots a_m}_{b_1\cdots b_m}\right)\!(u) \otimes 1\right)
\left(1 \otimes T\!\left(^{a_1\cdots a_m}_{b_1\cdots
b_m}\right)\!(v)\right)
\nonumber\\
&&=\left(1 \otimes T\!\left(^{a_1\cdots a_m}_{b_1\cdots b_m}\right)\!(v)\right)
\left(T\!\left(^{a_1\cdots a_m}_{b_1\cdots b_m}\right)\!(u) \otimes 1\right)
R^{(m)}_{12}(u-v)\;,
\end{eqnarray}
where $ T\!\left(^{a_1\cdots a_m}_{b_1\cdots
b_m}\right)\!(u)=\sum_{i,j=1}^m T_{a_ib_j}(u)\;E_{ij}$\;.
\begin{definition}
\label{qminor}
The quantum minor $t\!\left(^{a_1\cdots a_m}_{b_1\cdots
b_m}\right)\!(u)$ of $T(u)$ is defined by
\begin{eqnarray}
t\!\left(^{a_1\cdots a_m}_{b_1\cdots b_m}\right)\!(u)
&=&
qdet T\!\left(^{a_1\cdots a_m}_{b_1\cdots
b_m}\right)\!(u)
~=~\sum_{\sigma\in\fS_m}sgn(\sigma)~T_{a_{\sigma(1)},b_1}(u)\cdots
T_{a_{\sigma(m)},b_m}(u+m-1)
\end{eqnarray}
\end{definition}
One can show that
\begin{eqnarray}
t\!\left(^{a_1\cdots a_m}_{b_1\cdots
b_m}\right)\!(u)&=&
\sum_{\sigma\in\fS_m}
sgn(\sigma)~T_{a_{1},b_{\sigma(1)}}(u+m-1)\cdots
T_{a_{m},b_{\sigma(m)}}(u).
\end{eqnarray}  
By convention, when $m<1$, the quantum minor is equal to one.
Quantum minors satisfy some properties \cite{Molev2001} which are
analogous to those of numerical matrices
minors.
\begin{proposition}
The quantum minor $t\!\left(^{a_1\cdots a_m}_{b_1\cdots
b_m}\right)\!(u)$ verifies the following properties:
\begin{enumerate}
\item It is antisymmetric, i.e. for $\rho\in\fS_m$, 
\begin{eqnarray}
t\!\left(^{a_1\cdots a_m}_{b_1\cdots b_m}\right)\!(u)&=& 
sgn(\rho)~t\!\left(^{\rho(a_1)\cdots \rho(a_m)}_{b_1~~~\cdots b_m}\right)\!(u)
~=~
sgn(\rho)~t\!\left(^{a_1~~~\cdots a_m}_{\rho(b_1)\cdots \rho(b_m)}\right)\!(u).
\end{eqnarray}
\item It is alternated, i.e. if there exists
$i\neq j$ such that $a_i=a_j$ or $b_i=b_j$, 
then $t\!\left(^{a_1\cdots a_m}_{b_1\;\cdots
b_m}\right)\!(u)=0$.
\item It can be expanded with
respect to its last column or its last row as follows:
\begin{eqnarray}
\label{dev_col}
t\!\left(^{a_1\cdots a_m}_{b_1\cdots
b_m}\right)\!(u)&=&
\sum_{k=1}^{m}(-1)^{k+m}~
t\!\left(^{a_1\,\cdots\, a_{k-1}\,a_{k+1}\,\cdots
\,a_{m-1}\,a_{m}}_{b_1\;\cdots\, b_{k-1}\,b_{k}\,\hspace{3mm}\,\cdots\,
\,b_{m-2}\,b_{m-1}}\right)\!(u)
\;T_{a_k,b_m}(u+m-1)\\
\label{dev_row}
&=&\sum_{k=1}^{m}(-1)^{k+m}~
T_{a_m,b_k}(u+m-1)\;
t\!\left(^{a_1\,\cdots\, a_{k-1}\,a_{k}\,\hspace{3mm}\cdots
\,a_{m-2}\,a_{m-1}}_{b_1\;\cdots\, b_{k-1}\,b_{k+1}\,\,\cdots\,
\,b_{m-1}\,b_{m}}\right)\!(u).
\end{eqnarray}
\end{enumerate}
\end{proposition}
From the defining relations (\ref{RTT_relcom}), the commutation
relations of the quantum minors with $T_{i,j}(u)$ can be computed:
\begin{eqnarray}
\label{relcom_minor}
&&(u-v)\left[T_{i,j}(u),t\!\left(^{a_1\cdots a_m}_{b_1\cdots
b_m}\right)\!(v)\right]\nonumber\\
&&\hspace{10mm}=\sum_{k=1}^m \Big(
t\!\left(^{a_1\cdots a_{k-1}\,a_{k}a_{k+1}\cdots a_m}_
{b_1\cdots\, b_{k-1}\,j\hspace{2mm}\,b_{k+1}\cdots b_m}\right)\!(v)\;
T_{i,b_{k}}(u)
~-~
T_{a_{k},j}(u)\;t\!\left(^{a_1\cdots a_{k-1}\,i\hspace{2mm}a_{k+1}\cdots a_m}_
{b_1\cdots\, b_{k-1}\,b_k\,b_{k+1}\cdots b_m}\right)\!(v)\Big).
\end{eqnarray}
A corollary of (\ref{relcom_minor}) is that the quantum minor 
$t\!\left(^{a_1\cdots a_m}_{b_1\cdots b_m}\right)\!(u)$ 
lies in the centre of the subalgebra generated by
$\{T_{a_i,b_j}(u)|1\le i,j \le m\}$ i.e. 
\begin{eqnarray}
\label{comm}
\left[T_{a_i,b_j}(u)\;,\;
t\!\left(^{a_1\cdots a_m}_{b_1\cdots
b_m}\right)\!(v)\right]=0\;,\quad\mbox{for}\quad 1\le i,j \le m.
\end{eqnarray} 
The map 
\begin{eqnarray}
\label{tilde}
 T_{i,j}(u)\longmapsto t\!\left(^{1\cdots p\,p+i}
_{1\cdots p\,p+j}\right)\!(u)
\end{eqnarray}
defines an algebra homomorphism
$Y(sl_{n-p})\longrightarrow Y(sl_n)$, for $1\le p \le n-1$ and $1\le
i,j \le n-p$. Note that this homomorphism allows us to compute a simple
way the commutation relations among the $t\!\left(^{1\cdots p\,p+i}
_{1\cdots p\,p+j}\right)\!(u)$ minors.\\
Finally, quantum minors allow us to express some elements of the inverse
matrix of $T(u)$ thanks to the following proposition, proved by A.I.Molev \cite{Molev2001}:
\begin{proposition}
For $1\le i,j \le n$, the following equality holds
\begin{equation}
\Big(T^{-1}(u+n-1)\Big)_{i,j}=(-1)^{i+j}~t\!\left(^{1\cdots
j-1\;j+1\cdots n}_{1\cdots i-1\;\;i+1\cdots n}\right)\!(u).
\end{equation}
\end{proposition}

\section{Isomorphisms between the two realisations of $Y(sl_n)$}
\setcounter{equation}{0}

For clarity purposes, the isomorphism between the two previous realisations 
is recalled, see e.g. \cite{Dri88,Iohara}. Two presentations of this
isomorphism are possible. The first one uses the quantum minors and the
second one uses the Gauss decomposition.
\begin{theorem}
\label{iso}
The map~~$\phi:\cA\rightarrow U(R)$
\begin{eqnarray}
e_i\left(u+\frac{i-2}{2}\right)
&\mapsto&
\left(t\!\left(^{1\cdots i}_{1\cdots i}\right)\!(u)\right)^{-1}
t\!\left(^{1\cdots i-1\;i}_{1\cdots i-1\;i+1}\right)\!(u)\\
f_i\left(u+\frac{i-2}{2}\right)
&\mapsto&
t\!\left(^{1\cdots i-1\;i+1}_{1\cdots i-1\;i}\right)\!(u)
\left(t\!\left(^{1\cdots i}_{1\cdots i}\right)\!(u)\right)^{-1}\\
\label{phih1}
h_i\left(u+\frac{i-2}{2}\right)
&\mapsto&
\left(t\!\left(^{1\cdots i}_{1\cdots i}\right)\!(u)\right)^{-1}\;
t\!\left(^{1\cdots i-1}_{1\cdots i-1}\right)\!(u)\;
t\!\left(^{1\cdots i+1}_{1\cdots i+1}\right)\!(u-1)\;
\left(t\!\left(^{1\cdots i}_{1\cdots i}\right)\!(u-1)\right)^{-1}
\end{eqnarray}
is an algebra isomorphism.
\end{theorem} 
Note that the image of $h_i(u)$ can
be written differently  as:
\begin{eqnarray}
\label{phih2}
\phi\left(h_i\left(u+\frac{i-2}{2}\right)\right)\hspace{24mm}
&&\hspace{-30mm}
=
\left(t\!\left(^{1\cdots i}_{1\cdots i}\right)\!(u)\right)^{-1}\;
t\!\left(^{1\cdots i-1\;i+1}_{1\cdots i-1\;i+1}\right)\!(u)
-\phi\left(
f_i\left(u+\frac{i}{2}\right)
e_i\left(u+\frac{i-2}{2}\right)
\right)\\
\label{phih3}
&&\hspace{-30mm}=
t\!\left(^{1\cdots i-1\;i+1}_{1\cdots i-1\;i+1}\right)\!(u)\;
\left(t\!\left(^{1\cdots i}_{1\cdots i}\right)\!(u)\right)^{-1}
-\phi\left(
f_i\left(u+\frac{i-2}{2}\right)
e_i\left(u+\frac{i}{2}\right)
\right)
\end{eqnarray}
The other presentation of the isomorphism $\phi$
uses the Gauss decompositions of
the matrix $T(u)$:
\begin{eqnarray}
&&
\label{gauss1}
\hspace{-15mm}
T(u)=\left(\begin{array}{c c c c}
      1 &   &   & 0   \\
      f_{2,1}(u) & 1 &  &    \\
      \vdots  &\ddots   & \ddots &     \\
      f_{n,1}(u)  &   \cdots  & f_{n,n-1}(u)  & 1  
    \end{array}
  \right)
\!\!
\left(\begin{array}{c c c  }
      {k}_{1}(u) &   &     0   \\
       &\ddots   &         \\
      0 &        &  {k}_{n}(u)  
    \end{array}
  \right)
\!\!
\left(\begin{array}{c c c c}
    1 & {e}_{1,2}(u) 
& \cdots  & {e}_{1,n}(u) \\
    & 1 &\ddots   & \vdots \\ 
    & &\ddots   & {e}_{n-1,n}(u)  \\
    0  &   &   & 1  
    \end{array}
\right)
\\
\label{gauss2}
&&\hspace{-6mm}=
\left(\begin{array}{c c c c}
      1 & \widetilde{e}_{1,2}(u)     & \cdots  & \widetilde{e}_{1,n}(u) \\
        & 1 &\ddots   & \vdots \\ 
          & &\ddots   & \widetilde{e}_{n-1,n}(u)  \\
        0  &   &   & 1  
    \end{array}
\right)\!\!
\left(\begin{array}{c c c  }
      \widetilde{k}_{1}(u) &   &     0   \\
       &\ddots   &         \\
      0 &        &  \widetilde{k}_{n}(u)  
    \end{array}
  \right)\!\!
\left(\begin{array}{c c c c}
      1 &   &   & 0   \\
      \widetilde{f}_{2,1}(u) & 1 &  &    \\
      \vdots  &\ddots   & \ddots &     \\
      \widetilde{f}_{n,1}(u)  &   \cdots  & \widetilde{f}_{n,n-1}(u)  & 1  
    \end{array}
  \right)
\end{eqnarray}
The expression of the elements of the Gauss decomposition
(\ref{gauss1}) in terms of quantum minors has been
computed by K. Iohara \cite{Iohara}. For the alternative Gauss
decomposition (\ref{gauss2}), the computations are similar and one
obtains:
\begin{proposition}
\label{Gauss}
Let $1\le i < j \le n$ and $1\le p \le n$. The formal series
$e_{i,j}(u),\;
\widetilde{e}_{i,j}(u),\;f_{j,i}(u),\;\widetilde{f}_{j,i}(u),\;
$ $k_p(u)$ and $\widetilde{k}_p(u)$ in the Gauss
decompositions, can be expressed in terms of quantum minors:
\begin{eqnarray}
&&
{e}_{i,j}(u+i-1)=
\left(t\!\left(^{1\cdots i}_{1\cdots i}\right)\!(u)\right)^{-1}
t\!\left(^{1\cdots i-1\;i}_{1\cdots i-1\;j}\right)\!(u)\;,\\
&&
{f}_{j,i}(u+i-1)=
t\!\left(^{1\cdots i-1\;j}_{1\cdots i-1\;i}\right)\!(u)
\left(t\!\left(^{1\cdots i}_{1\cdots i}\right)\!(u)\right)^{-1}\;,\\
&&
{k}_j(u+j-1)=
t\!\left(^{1\cdots j}_{1\cdots j}\right)\!(u)
\left(t\!\left(^{1\cdots j-1}_{1\cdots j-1}\right)\!(u)\right)^{-1}\;,\\
&&
{k}_1(u)=t\!\left(^{1}_{1}\right)\!(u)=T_{1,1}(u)\;,
\end{eqnarray}
and
\begin{eqnarray}
&&\widetilde{e}_{i,j}(u+n-j)=
t\!\left(^{i\;j+1\cdots n}_{j\;j+1\cdots n}\right)\!(u)
\left(t\!\left(^{j+1\cdots n}_{j+1\cdots n}\right)\!(u)\right)^{-1}\;,\\
&&
\widetilde{f}_{j,i}(u+n-j)=
\left(t\!\left(^{j+1\cdots n}_{j+1\cdots n}\right)\!(u)\right)^{-1}
t\!\left(^{j\;j+1\cdots n}_{i\;j+1\cdots n}\right)\!(u)\;,\\
&&\widetilde{k}_i(u+n-i)=
\left(t\!\left(^{i+1\cdots n}_{i+1\cdots n}\right)\!(u)\right)^{-1}
t\!\left(^{i\cdots n}_{i\cdots n}\right)\!(u)\;,
\\
&&\widetilde{k}_n(u)=t\!\left(^{n}_{n}\right)\!(u)=T_{n,n}(u)\;.
\end{eqnarray}
\end{proposition}
{\it Remark:} Proposition \ref{Gauss} proves the
existence of the two Gauss decomposition.\\
Then, proposition \ref{Gauss} implies that the 
map~~$\widetilde{\phi}:\cA \longrightarrow U(R)$
\begin{eqnarray}
e_i(u)&\longmapsto&{e}_{i,i+1}\left(u+\frac{i}{2}\right)
\\
f_i(u)&\longmapsto&{f}_{i+1,i}\left(u+\frac{i}{2}\right)
\\
h_i(u)&\longmapsto&{k}_{i+1}\left(u+\frac{i}{2}\right)
{k}_{i}^{-1}\left(u+\frac{i}{2}\right)
\end{eqnarray}
is an algebra isomorphism, for $1\le i \le n-1$.\\
In the following, as in equation (\ref{inverse}), $T^*(u)$ denotes
$T(-u)^{-1}$ and $t^*\!\left(^{a_1\cdots a_m}_
{b_1\cdots b_m}\right)\!(u)$ denotes the quantum minors of
$T^*(u)$ .
\begin{corollary}
\label{inverse_minor}
For $1\le m \le n$, the following equalities hold
\begin{eqnarray}
t^*\!\left(^{1\cdots m}_{1\cdots m}\right)\!(-u-n+1)
&=& 
t\!\left(^{m+1\cdots n}_{m+1\cdots n}\right)\!(u)\\
t^*\!\left(^{1\cdots m-1\;m+1}_{1\cdots m-1\;m}\right)\!(-u-n+1)
&=& -
t\!\left(^{m+1\;m+2\cdots n}_{m~~~m+2\cdots n}\right)\!(u)\\
t^*\!\left(^{1\cdots m-1\;m}_{1\cdots m-1\;m+1}\right)\!(-u-n+1)
&=& -
t\!\left(^{m~~~m+2\cdots n}_{m+1\;m+2\cdots n}\right)\!(u)\\
t^*\!\left(^{1\cdots m-1\;m+1}_{1\cdots m-1\;m+1}\right)\!(-u-n+1)
&=& 
t\!\left(^{m\;m+2\cdots n}_{m\;m+2\cdots n}\right)\!(u)
\end{eqnarray}
\end{corollary}
\textbf{Proof:} Let $T(u)$ decomposed according to 
(\ref{gauss2}). Then, $T^*(u)$ decomposes as in
(\ref{gauss1}). Using the relation $T^*(u)=T(-u)^{-1}$, we deduce 
\begin{eqnarray}
e^*_{i,i+1}(u)&=&-e_{i,i+1}(-u)\;,\\
f^*_{i+1,i}(u)&=&-f_{i+1,i}(-u)\;,\quad\mbox{for}\quad 1\le i\le n-1\\
k^*_{i}(u)&=&(k_{i}(-u))^{-1}\;,\quad\mbox{for}\quad 1\le i\le n
\end{eqnarray}
with obvious notations. Finally, using proposition \ref{Gauss},
the equalities are proven.
\finproof\\
For $1 \le i < j \le n$, the elements $e_{i,j}^{(0)}=T_{i,j}^{(1)}$
and $f_{j,i}^{(0)}=T_{j,i}^{(1)}$ are root generators of the algebra $sl_n$,
which can be expressed in terms of simple root generators, $e_i^{(0)}$ and
$f_i^{(0)}$, as follows:
\begin{eqnarray}
\label{roote}
e_{i,j}^{(0)}&=&[\cdots[e_{j-1}^{(0)}\;,\;e_{j-2}^{(0)}],\;e_{j-3}^{(0)}],\cdots
],e_{i+1}^{(0)}],e_i^{(0)}]\;,\\
\label{rootf}
f_{j,i}^{(0)}&=&[f_i^{(0)},[f_{i+1}^{(0)}\;,[\cdots,[f_{j-3}^{(0)}\;,
[f_{j-2}^{(0)}\;,\;f_{j-1}^{(0)}]\cdots].
\end{eqnarray}
{\it Remark:} In (\ref{roote}) and (\ref{rootf}), the isomorphism
$\widetilde{\phi}$ has been omitted for simplicity. In the
following, this losely notation is always used, i.e. the
isomorphisms between two realisations of the same algebra are omitted.

\section{The Hopf structure of $Y(sl_n)$ in the Drinfel'd basis}
\setcounter{equation}{0}

Before giving the Hopf structure of $Y(sl_n)$ in the Drinfel'd basis,
the images of any quantum minor under the coproduct, the antipode and the
counit are needed.
Let us recall that $T^*(u)$ denotes
$T(-u)^{-1}$ and $t^*\!\left(^{a_1\cdots a_m}_
{b_1\cdots b_m}\right)\!(u)$ denotes its quantum minors.
\begin{proposition}
Let $1\le m\le n$, $1\le a_1< \cdots < a_m \le n$ and $1\le b_1<
\cdots < b_m \le n$.
The images of a quantum minor under the coproduct, the antipode and the
counit are given by:
\begin{eqnarray}
&&\Delta\left(t\!\left(^{a_1\cdots a_m}_
{b_1\cdots b_m}\right)\!(u)\right)
= \sum_{1\le c_1<\cdots< c_m \le n}
t\!\left(^{a_1\cdots a_m}_{c_1\cdots c_m}\right)\!(u)
\otimes 
t\!\left(^{c_1\cdots c_m}_{b_1\cdots b_m}\right)\!(u)\;,\\
&&S\left(t\!\left(^{a_1\cdots a_m}_
{b_1\cdots b_m}\right)\!(u)\right)
=(-1)^{\left[\frac{i}{2}\right]}~t^*\!\left(^{a_1\cdots a_m}_
{b_1\cdots b_m}\right)\!(-u-m+1)\;,\\
&&\eps\left(t\!\left(^{a_1\cdots a_m}_
{b_1\cdots b_m}\right)\!(u)\right)
=\delta_{a_1,b_1}\cdots \delta_{a_m,b_m}\;,
\end{eqnarray}
where $\left[\frac{i}{2}\right]$ is the integer part of $\frac{i}{2}$.
\end{proposition}
\textbf{Proof:} Direct computation, using the definition of the
Hopf structure (\ref{hopf}) and the property that the
comultiplication and counit are morphisms while the antipode is
an anti-morphism.
\finproof\\
In particular, one obtains the following well-known result:
\begin{eqnarray}
\Delta(qdetT(u))=qdetT(u)\;\otimes\;qdetT(u).
\end{eqnarray}
The comultiplication, the antipode and the counit are established in
the Drinfel'd basis thanks to the
isomorphism $\phi$ (see theorem \ref{iso}) between $\cA$ and
$U(R)$.

\subsection{Comultiplication}
The adjoint actions of the elements of the algebra $sl_n$ on $X\in Y(sl_n)$
will be denoted by, for $1 \le i< j \le n$:
\begin{eqnarray}
\label{ad}
Ad^\pm_{e_{i,j}}(X)=\pm[e_{i,j}^{(0)}\;,X]\quad \mbox{and}\quad
Ad^\pm_{f_{j,i}}(X)=\pm[f_{j,i}^{(0)}\;,X]\;.
\end{eqnarray}
Moreover, by convention $Ad^\pm_{e_{i,i}}(X)=X$ and
$Ad^\pm_{f_{i,i}}(X)=X$.
To determine the explicit form of the comultiplication,
the following generalisation
of the adjoint action, depending on a spectral parameter, is useful.
\begin{definition}
Let $1 \le i\le j \le n$,  $1 \le \alpha \le
n$ and $X$ an element of $Y(sl_n)$.
The generalised adjoint actions are defined by
\begin{eqnarray}
\label{actione1}
{^\alpha}\cE_{i,j}(u)(X)&=&
Ad^+_{e_{i,j}}(X)
+\delta_{i\le\alpha<j}~
Ad^-_{e_{i,\alpha}}\left(Ad^+_{e_{\alpha+1,j}}(e_{\alpha}(u))\right)~
X
\;,
\\
\label{actione2}
\cE^\alpha_{i,j}(u)(X)&=&
Ad^+_{e_{i,j}}(X)
+\delta_{i\le\alpha<j}~X~
Ad^-_{e_{i,\alpha}}\left(Ad^+_{e_{\alpha+1,j}}(e_{\alpha}(u))\right)
\;,
\end{eqnarray}
and
\begin{eqnarray}
\label{actionf1}
{^\alpha}\cF_{j,i}(u)(X)&=&
Ad^-_{f_{j,i}}(X)
+\delta_{i\le\alpha<j}~
Ad^+_{f_{\alpha,i}}\left(Ad^-_{f_{j,\alpha+1}}(f_{\alpha}(u))\right)~X\;,
\\
\cF^\alpha_{j,i}(u)(X)&=&
Ad^-_{f_{j,i}}(X)
+\delta_{i\le\alpha<j}~
X~
Ad^+_{f_{\alpha,i}}\left(Ad^-_{f_{j,\alpha+1}}(f_{\alpha}(u))\right)\;,
\end{eqnarray}
where~~$\delta_{i\le\alpha<j-1}=\begin{cases} 1
&\mbox{if}~~i\le\alpha<j-1\;,\\ 0 & otherwise\;. \end{cases}$
\end{definition}
Let ${^\alpha}\cG$, $\cG^\beta$ be any actions on
$Y(sl_n)$. Hereafter, for
simplicity, the notation ${^\alpha}\cG^\beta$ means either
${^\alpha}\cG$ or $\cG^\beta$.\\
To compute 
the comultiplication, we also need:
\begin{definition}
For $1 \le m\le n$, $1 \le k_1 < k_2 <
\cdots < k_m \le n$ and $k_m\neq m$, ${^\alpha}E^\beta_{
k_1,k_2,\cdots,k_m}(u)$ and ${^\alpha}F^\beta_{
k_1,k_2,\cdots,k_m}(u)(X)$ are defined by for $X\in Y(sl_n)$:
\begin{eqnarray}
{^\alpha}E^\beta_{
k_1,k_2,\cdots,k_m}(u)(X)&=&
\left(
\prod_{1\le i \le m-1}^{\longrightarrow}
{^\alpha}\cE^\beta_{i,k_i}(u)
\right)
\left(
{^\alpha}\cE^\beta_{m+1,k_m}(u)\;
\left(X\right)\right)\;,\\
{^\alpha}F^\beta_{
k_1,k_2,\cdots,k_m}(u)(X)&=&
\left(
\prod_{1\le i \le m-1}^{\longrightarrow}
{^\alpha}\cF^\beta_{k_i,i}(u)
\right)
\left(
{^\alpha}\cF^\beta_{k_m,m+1}(u)\;
\left(X\right)\right)\;,
\end{eqnarray}
where, for $\{\cG_p|1\le p \le m-1\}$, a set of actions on $Y(sl_n)$,
we denote
\begin{eqnarray}
\prod^{\longrightarrow}_{1\le i \le m-1}\cG_i~(X)
=\cG_1(\cdots(\cG_{m-2}(\cG_{m-1}~(X))\cdots)\;.
\end{eqnarray}
\end{definition}
In particular, one gets for $k>1$
\begin{eqnarray}
{^\alpha}E^\beta_k(u)(X)=
{^\alpha}\cE^\beta_{2,k}(u)(X)
&,&
{^\alpha}F^\beta_k(u)(X)={^\alpha}\cF^\beta_{k,2}(u)(X).
\end{eqnarray}
By convention, if the set of indices 
$\{k_1,k_2,\cdots ,k_m\}$ is empty, then
${^\alpha}E^\beta_{
k_1,k_2,\cdots,k_m}(u)(X)=1$ and ${^\alpha}F^\beta_{
k_1,k_2,\cdots,k_m}(u)(X)=1$. Remark that these generators can be
expressed only in terms of the elements of the Drinfel'd basis, thanks
to equations (\ref{roote}) and (\ref{rootf}).\\
These generalised actions show up in the following lemma:
\begin{lemma}
For $1\le i \le n-1$,
$1\le a_1 <\; \cdots \;< a_i \le n$ and $a_i \neq i$, one gets
\label{lem_1}
\begin{eqnarray}
\label{iXp}
\Big(t\!\left(^{1\cdots i}_{1\cdots i}\right)\!(u)\Big)^{-1}
t\!\left(^{1~\cdots \;i}_{a_1\cdots \;a_{i}}\right)\!(u)
&=&
{^i}E_{a_1,\cdots,a_{i}}\left(u+\frac{i-2}{2}\right)
\left(e_i\left(u+\frac{i-2}{2}\right)\right)\\
\label{iXpb}
t\!\left(^{1~\cdots \;i}_{a_1\cdots \;a_{i}}\right)\!(u)
\Big(t\!\left(^{1\cdots i}_{1\cdots i}\right)\!(u)\Big)^{-1}
&=&
E^i_{a_1,\cdots,a_{i}}\left(u+\frac{i}{2}\right)
\left(e_i\left(u+\frac{i}{2}\right)\right)\\
\label{iXm}
~\Big(t\!\left(^{1\cdots i}_{1\cdots i}\right)\!(u)\Big)^{-1}
t\!\left(^{a_1\cdots a_i}_{1~\cdots i}\right)\!(u)
&=&
{^i}F_{a_1,\cdots,a_{i}}\left(u+\frac{i}{2}\right)
\left(f_i\left(u+\frac{i}{2}\right)\right)\\
\label{iXmb}
t\!\left(^{a_1\cdots a_i}_{1~\cdots i}\right)\!(u)
~\Big(t\!\left(^{1\cdots i}_{1\cdots i}\right)\!(u)\Big)^{-1}
&=&
F^i_{a_1,\cdots,a_{i}}\left(u+\frac{i-2}{2}\right)
\left(f_i\left(u+\frac{i-2}{2}\right)\right)\\
\label{iHp}
\Big(t\!\left(^{1\cdots i}_{1\cdots i}\right)\!(u)\Big)^{-1}
t\!\left(^{1~\cdots i-1~~i+1}_{a_1\cdots a_{i-1}\;a_{i}}\right)\!(u)
&=&
{^i}E_{a_1,\cdots,a_{i}}\left(u+\frac{i-2}{2}\right)
\left(\widetilde{g}{_i}\left(u+\frac{i-2}{2}\right)
\right)\\
\label{iHpb}
t\!\left(^{1~\cdots i-1~~i+1}_{a_1\cdots a_{i-1}\;a_{i}}\right)\!(u)
\Big(t\!\left(^{1\cdots i}_{1\cdots i}\right)\!(u)\Big)^{-1}
&=&
E^i_{a_1,\cdots,a_{i}}\left(u+\frac{i}{2}\right)
\left(g_i\left(u+\frac{i-2}{2}\right)\right)\\
\label{iHm}
~\Big(t\!\left(^{1\cdots i}_{1\cdots i}\right)\!(u)\Big)^{-1}
t\!\left(^{a_1\cdots a_{i-1}\;a_i}_{1~\cdots i-1~~i+1}\right)\!(u)
&=&
{^i}F_{a_1,\cdots,a_{i}}\left(u+\frac{i}{2}\right)
\left(\widetilde{g}_i\left(u+\frac{i-2}{2}\right)\right)\\
\label{iHmb}
t\!\left(^{a_1\cdots a_{i-1}\;a_i}_{1~\cdots i-1~~i+1}\right)\!(u)
~\Big(t\!\left(^{1\cdots i}_{1\cdots i}\right)\!(u)\Big)^{-1}
&=&
F^i_{a_1,\cdots,a_{i}}\left(u+\frac{i-2}{2}\right)
\left(g_i\left(u+\frac{i-2}{2}\right)\right)
\end{eqnarray}
where
\begin{eqnarray}
g_i(u)&=&h_i\left(u\right)+
f_i\left(u\right)
e_i\left(u+1\right)\\
\widetilde{g}_i(u)&=&
h_i\left(u\right)+
f_i\left(u+1\right)
e_i\left(u\right)
\end{eqnarray}
\end{lemma}
\textbf{Proof:}
The proof is only given for (\ref{iXp}).
Let $i$ and $\{a_1, \cdots ,a_i\}$ fixed as in the lemma.
The first step is to evaluate the quantum minor
$t\!\left(^{1~\cdots i}_{a_1 \cdots a_i}\right)\!(u)$ in terms of
the quantum minor 
$t\!\left(^{1\cdots i-1\;i}_{1 \cdots i-1\;i+1}\right)\!(u)$ and 
in terms of some generators of the $sl_{n}$ algebra.
Selecting the coefficient of $u^0v^{-1}$ in equation
(\ref{relcom_minor}), the following relation is
obtained for $1\le b_1 <\cdots<b_i\le n$, $1\le p\le i$ and $1\le m\le
n-b_p$:
\begin{eqnarray}
Ad^+_{e_{b_p,b_p+m}}
\left(t\!\left(^{1~\cdots p~ \cdots i}_{b_1 \cdots b_p
\cdots b_i}\right)\!(v)\right)
=
\left[T_{b_p,b_p+m}^{(1)}\;,\;
t\!\left(^{1~\cdots p~ \cdots i}_{b_1 \cdots b_p
\cdots b_i}\right)\!(v)\right]
=t\!\left(^{1~\cdots p~~~~ \cdots i}
_{b_1 \cdots b_p+m \cdots b_i}\right)\!(v).
\end{eqnarray} 
This relation allows us to increase the parameters of the
quantum minor. 
Using this relation, the indices $\{1,\cdots,i-1,\;i+1\}$
of 
$t\!\left(^{1\cdots i-1\;i}_{1 \cdots i-1\;i+1}\right)\!(u)$ 
can be increased up to
$\{a_1,\cdots,a_i\}$:
\begin{eqnarray}
\label{A1}
\left(\prod_{1\le p \le i-1}^{\longrightarrow}
Ad^+_{e_{p,a_p}}
\right)\left(
Ad^+_{e_{i+1,a_i}}
\left(t\left(^{1\cdots i-1\;i}_{1 \cdots i-1\;i+1}\right)\!(u)\right)
\right)
=
t\!\left(^{1~\cdots i}_{a_1 \cdots a_i}\right)\!(u)
\end{eqnarray} 
The second step of the proof consists in determining the commutator of 
$\Big(t\!\left(^{1\cdots i}_{1\cdots i}\right)\!(u)\Big)^{-1}$ with 
$e_{j,j+1}^{(0)}$. The commutation 
relations are computed using equation (\ref{relcom_minor}):
\begin{eqnarray}
\label{rel_com_modu}
(u-v) 
[T_{j,j+1}(u)\;,\;t\!\left(^{1\cdots i}_{1\cdots i}\right)\!(v)]
=
\delta_{ij}~t\!\left(^{1\cdots i-1\;i}
_{1\cdots i-1\;i+1}\right)\!(v)T_{ii}(u)
+O\left(\frac{1}{u}\right)
\end{eqnarray}
Multiplying by 
$\Big(t\!\left(^{1\cdots i}_{1\cdots i}\right)\!(v)\Big)^{-1}$ 
both sides of 
the $u^0$ coefficient in (\ref{rel_com_modu}), 
one gets
\begin{eqnarray}
\left[e_{j,j+1}^{(0)}\;,\;
\Big(t\!\left(^{1\cdots i}_{1\cdots i}\right)\!(v)\Big)^{-1}\right]
&=&
\label{com2}
-\delta_{ij}~
e_i\left(v+\frac{i-2}{2}\right)
\Big(t\!\left(^{1\cdots i}_{1\cdots i}\right)\!(v)\Big)^{-1}
\end{eqnarray}
Thus, thanks to the relations (\ref{roote}) and (\ref{com2}), one obtains
for $X\in Y(sl_n)$:
\begin{eqnarray}
\Big(t\!\left(^{1\cdots i}_{1\cdots i}\right)\!(u)\Big)^{-1}
Ad^+_{e_{p,a_p}}(X)
={^i}\cE_{p,a_p}\left(
u+\frac{i-2}{2}
\right)
\left(
\Big(t\!\left(^{1\cdots i}_{1\cdots i}\right)\!(u)\Big)^{-1}X
\right)
\end{eqnarray}
This proves the equation (\ref{iXp}).
Equation (\ref{iXpb}) is proven along the same lines, remarking that
\begin{eqnarray}
t\!\left(^{1\cdots i-1\;i}_{1\cdots i-1\;i+1}\right)\!(u-1)
\Big(t\!\left(^{1\cdots i}_{1\cdots i}\right)\!(u-1)\Big)^{-1}
&=&
\Big(t\!\left(^{1\cdots i}_{1\cdots i}\right)\!(u)\Big)^{-1}
t\!\left(^{1\cdots i-1\;i}_{1\cdots i-1\;i+1}\right)\!(u),
\end{eqnarray}
which explains the shift in the spectral parameter between (\ref{iXp}) 
and (\ref{iXpb}). 
For the relations (\ref{iXm})-(\ref{iHmb}), the proof is analogous.
\finproof\\
Now we can state the main theorem of the letter.
\begin{theorem}
\label{coproduct}Let $1\le i \le n-1$.
The comultiplication in the Drinfel'd basis is given by: 
\begin{eqnarray}
\label{copx+}
&&\hspace{-5mm}\Delta(e_i(u))=
\sum_{m=0}^{+\infty}(-1)^{m}
\left(
\sum_{\substack{1\le b_1<\cdots<b_i \le n\\b_i \neq i}}
~~{^i}E_{b_1,\cdots,b_i}(u)(e_i(u))
\otimes
{^i}F_{b_1,\cdots, b_i}(u+1)(f_i(u+1))
\right)^m\nonumber\\
&&\hspace{17mm}\times\left(
1\otimes e_i(u)+
\sum_{\substack{1\le a_1<\cdots<a_i\le n\\a_i\neq i}}
{^i}E_{a_1, \cdots, a_i}(u)(e_i(u))~
\otimes~
{^i}F_{a_1,\cdots, a_i}(u+1)\left(\widetilde{g}_i(u)\right)
\right)
\\
&&\hspace{-5mm}\Delta(f_i(u))=
\left(
f_i(u)\otimes 1+
\sum_{\substack{1\le a_1<\cdots<a_i \le n\\a_i\neq i}}
E^i_{a_1,\cdots,a_i}(u+1)(g_i(u))
\otimes F^i(a_1 \cdots a_i)(u)(f_i(u))
\right)
\nonumber\\
&&
\hspace{17mm}\times\sum_{m=0}^{+\infty}(-1)^{m}
\left(
\sum_{\substack{1\le b_1<\cdots<b_i\le n\\b_i \neq i}}
~~E^i_{b_1,\cdots,b_i}(u+1)(e_i(u+1))
~~\otimes~~
F^i_{b_1,\cdots,b_i}(u)(f_i(u))
\right)^m\\
&&\hspace{-5mm}\Delta(h_i(u))=
\left(
f_i(u)\otimes e_i(u+1)+\!\!
\sum_{\substack{1\le a_1<\cdots<a_i\le n\\a_i \neq i}}
E^i_{a_1,\cdots, a_i}(u+1)(g_i(u))
\otimes
F^i_{a_1,\cdots,a_i}(u)(g_i(u))
\right)
\nonumber\\
&& 
\hspace{17mm}\times\sum_{m=0}^{+\infty}(-1)^m
\left(
\sum_{\substack{1\le b_1<\cdots<b_i\le n\\b_i \neq i}}
~~E^i_{b_1,\cdots,b_i}(u+1)(e_i(u+1))
~~\otimes~~
F^i_{b_1,\cdots,b_i}(u+1)(f_i(u))
\right)^m\nonumber\\
&&\hspace{17mm}-\Delta(f_i(u))\Delta(e_i(u+1))
\end{eqnarray}
\end{theorem}
\textbf{Proof:}
The full proof is presented only for 
$e_i(u)$, the outline of proofs for $f_i(u)$ 
and $h_i(u)$ being similar.
The comultiplication in the Drinfel'd realisation 
is constructed thanks to
the isomorphism given in the theorem~\ref{iso}.
\begin{eqnarray}
\Delta\left(e_i\left(u+\frac{i-2}{2}\right)\right)&=&
\Delta\Big(t\!\left(^{1\cdots  i}_{1\cdots i}\right)\!(u)
\Big)^{-1}
\Delta\Big(t\!\left(^{1\cdots i-1\;i}_{1\cdots
i-1\;i+1}\right)\!(u)\Big)\\
&&\hspace{-50mm}=
\left(\sum_{b_1<\cdots<b_i} 
t\!\left(^{1~\cdots i}_{b_1\cdots b_i}\right)\!(u)
\otimes
t\!\left(^{b_1\cdots b_i}_{1~\cdots i}\right)\!(u)
\right)^{-1}
\sum_{a_1<\cdots<a_i} 
t\!\left(^{1~\cdots i}_{a_1\cdots
a_{i}}\right)\!(u)
\otimes
t\!\left(^{a_1\cdots a_{i-1}\;a_{i}}_{1~\cdots 
i-1~\;i+1}\right)\!(u)\\
&&\hspace{-50mm}=
\sum_{m\le 0}(-1)^m
\left(
\sum_{\substack{b_1<\cdots<b_i\\b_i \neq i}}
\Big(t\!\left(^{1\cdots i}_{1\cdots i}\right)\!(u)\Big)^{-1}
t\!\left(^{1~\cdots i}_{b_1\cdots b_i}\right)\!(u)
\otimes
\Big(t\!\left(^{1\cdots i}_{1\cdots i}\right)\!(u)\Big)^{-1}
t\!\left(^{b_1\cdots b_i}_{1~\cdots i}\right)\!(u)
\right)^m
\nonumber\\
&&\hspace{-50mm}\times
\label{cop_dem3}
\sum_{a_1<\cdots<a_i} 
\Big(t\!\left(^{1\cdots i}_{1\cdots i}\right)\!(u)\Big)^{-1}
t\!\left(^{1~\cdots \;i}_{a_1\cdots
\;a_{i}}\right)\!(u)
\otimes
\Big(t\!\left(^{1\cdots i}_{1\cdots i}\right)\!(u)\Big)^{-1}
t\!\left(^{a_1\cdots a_{i-1}\;a_{i}}_{1~\cdots
i-1~\;i+1}\right)\!(u)
\end{eqnarray}
The lemma \ref{lem_1} allows us to evaluate all the terms 
of equation (\ref{cop_dem3}), which proves (\ref{copx+}).
\finproof\\

\subsection{Antipode and counit}
As for the comultiplication, generalised adjoint actions must be
introduced to express the antipode.
\begin{definition}
Let $1 \le i,\;j \le n$,  $1 \le \alpha \le
n$. $X$ denotes an element of $Y(sl_n)$.
The generalised adjoint actions are defined by
\begin{eqnarray}
\label{actionh1}
{^\alpha}\cH_{i,j}(u)(X)&=&
\begin{cases}
1& \mbox{if}~~i>j\;\\
X& \mbox{if}~~i=j\;,\\
X~+~{^\alpha}\cE_{i,j}(u)
\Big({^\alpha}\cF_{j,i}(u+1)(X)\Big)
& otherwise,
\end{cases}\\
\label{actionh2}
\cH^\alpha_{i,j}(u)(X)&=&
\begin{cases}
1& \mbox{if}~~i>j\;\\
X& \mbox{if}~~i=j\;,\\
X~+~\cE^\alpha_{i,j}(u+1)
\Big(\cF^\alpha_{j,i}(u)(X)\Big)
& otherwise\;,
\end{cases}
\end{eqnarray} 
Let $1\le m\le n-1$ and $1\le k_1 <k_2<\cdots<k_m\le n$. Then, one has: 
\begin{eqnarray}
{^\alpha}H^\beta_{
k_1,k_2,\cdots,k_m}(u)(X)
&=&
\left(
\prod_{1\le p \le m-1}^{\longrightarrow}
{^\alpha}\cH^\beta_{p,k_p}(u)
\right)
\left(
{^\alpha}\cH^\beta_{m+1,k_m}(u)\;
\left(X\right)\right)\;.
\end{eqnarray}
Similary, $\widehat{E}_{m}(u)(X)$ is defined as:
\begin{eqnarray}
\begin{cases}
\cE_{2,n}^1(u+1)
\Big(\cF_{n-1,1}^1(u)(X)\Big)&\mbox{if}~~m=1\\
\cE_{1,n-m+1}^m(u+1)
\cF_{n-m,1}^m(u)\left(
\displaystyle\prod_{2\le p \le m-1}^{\longrightarrow}
\cH^m_{p,n-m+p}(u)
\right)\cE_{m+1,n}^m(u+1)
\cF_{n,m}^m(u)(X)
&otherwise,
\end{cases}
\end{eqnarray} 
and $\widehat{F}_{m}(u)(X)$ as:
\begin{eqnarray}
\begin{cases}
{^1}\cE_{1,n-1}(u)
\Big({^1}\cF_{n,2}(u+1)(X)\Big)&\mbox{if}~~m=1\\
{^m}\cE_{1,n-m}(u+1)
{^m}\cF_{n-m+1,1}(u)\left(
\displaystyle\prod_{2\le p \le m-1}^{\longrightarrow}
{^m}\cH_{p,n-m+p}(u)
\right)
{^m}\cE_{m,n}(u+1)
{^m}\cF_{n,m+1}(u)(X)&otherwise\;.
\end{cases}
\end{eqnarray} 
\end{definition}
To find the image under the antipode, the following lemma is needed.
\begin{lemma}
\label{lem_2}
For $1\le i \le n-1$ and $1 \le j \le i+1$, one has:
\begin{eqnarray}
t\!\left(^{j\;i+2\cdots n}_{j\;i+2\cdots n}\right)\!(u)
\Big(t\!\left(^{1\cdots n-i}_{1\cdots n-i}\right)\!(u)\Big)^{-1}
&=&
H^{n-i}_{j,i+2,\cdots,n}\left(u+\frac{n-i-2}{2}\right)\!\!
\left(
g_{n-i}\left(u+\frac{n-i-2}{2}\right)
\right)\label{lem_21}\\
\Big(t\!\left(^{1\cdots n-i}_{1\cdots n-i}\right)\!(u)\Big)^{-1}
t\!\left(^{j\;i+2\cdots n}_{j\;i+2\cdots n}\right)\!(u)
&=&
{^{n-i}}H_{j,i+2,\cdots,n}\left(u+\frac{n-i-2}{2}\right)\!\!
\left(
\widetilde{g}{_{n-i}}\left(u+\frac{n-i-2}{2}\right)
\right)\label{lem_22}
\end{eqnarray}
and
\begin{eqnarray}
t\!\left(^{i~~~i+2\cdots n}_{i+1\;i+2\cdots n}\right)\!(u)
\Big(t\!\left(^{1\cdots n-i}_{1\cdots n-i}\right)\!(u)\Big)^{-1}
&=&
\widehat{E}_{n-i}\left(u+\frac{n-i-2}{2}\right)
\left(e_{n-i}\left(u+\frac{n-i}{2}\right)\right)
\label{lem_23}\\
\Big(t\!\left(^{1\cdots n-i}_{1\cdots n-i}\right)\!(u)\Big)^{-1}
t\!\left(^{i+1\;i+2\cdots n}_{i~~~\;i+2\cdots n}\right)\!(u)
&=&
\widehat{F}_{n-i}\left(u+\frac{n-i-2}{2}\right)
\left(f_{n-i}\left(u+\frac{n-i}{2}\right)\right)
\label{lem_24}
\end{eqnarray}
\end{lemma}
\textbf{Proof:} This lemma is proven along the same lines as the lemma
\ref{lem_1}. 
\finproof\\
\begin{theorem}
\label{th_antipode}
The antipode and the counit in the Drinfel'd basis are given by, 
for $1\le i \le n-1$:
\begin{eqnarray}
S\left(e_{i}\left(u+\frac{n}{2}\right)\right)&=&-
\widehat{E}_{n-i}\left(u\right)
\left(
e_{n-i}\left(u+1\right)
\right)~~
\Big(
H_{i+1,\cdots,n}^{n-i}
\left(u\right)
\left(
g_{n-i}\left(u\right)
\right)
\Big)^{-1}
\\
S\left(f_{i}\left(u+\frac{n}{2}\right)\right)&=&-
\Big(
{^{n-i}}H_{i+1,i+2,\cdots,n}(u)
\left(
\widetilde{g}{_{n-i}}\left(u\right)
\right)\Big)^{-1}~~
\widehat{F}_{n-i}(u)
\left(
f_{n-i}\left(u+1\right)
\right)\\
S\left(h_{i}\left(u+\frac{n}{2}\right)\right)&=&
\Big(
{^{n-i}}H_{i+1,i+2,\cdots,n}(u)\left(\widetilde{g}_{n-i}(u)\right)
\Big)^{-1}~~
{^{n-i}}H_{i,i+2,\cdots,n}(u)\left(\widetilde{g}_{n-i}(u)\right)\nonumber\\
&&-~S(e_{n-i}(u+1))~S(f_{n-i}(u))
\end{eqnarray}
and
\begin{equation}
\epsilon(e_i(u))=0\;,\quad
\epsilon(f_i(u))=0\;,\quad
\epsilon(h_i(u))=1\;.
\end{equation}
\end{theorem}
\textbf{Proof:} The proof is given for
$S(e_{n-i}(u))$. $S(h_{n-i}(u))$
 and $S(f_{n-i}(u))$ are proven analogously.
\begin{eqnarray}
S\left(e_{i}\left(u+\frac{i-2}{2}\right)\right)
&=&
S\left(t\!\left(^{1\cdots i-1\;i}_{1\cdots i-1\;i+1}\right)\!(u)\right)
S\left(\left(t\!\left(^{1\cdots i}_{1\cdots
i}\right)\!(u)\right)\right)^{-1}\\
&=&
t^*\!\left(^{1\cdots i-1\;i}_{1\cdots i-1\;i+1}\right)\!(u)
\left(t^*\!\left(^{1\cdots i}_{1\cdots
i}\right)\!(u)\right)^{-1}\\
&=&-t\!\left(^{i~~~i+2\cdots n}_{i+1\;i+2\cdots n}\right)\!(u-n+i)
\Big(t\!\left(^{1\cdots n-i}_{1\cdots n-i}\right)\!(u-n+i)\Big)^{-1}\nonumber\\
&&\label{dem_S_3}
\hspace{-20mm}\times
\left(
t\!\left(^{i+1\cdots n}_{i+1\cdots n}\right)\!(u-n+i)
\Big(t\!\left(^{1\cdots n-i}_{1\cdots n-i}\right)\!(u-n+i)\Big)^{-1}
\right)^{-1}~~
(\mbox{corollary}(\ref{inverse_minor}))
\end{eqnarray}
The terms in the relation (\ref{dem_S_3}) can be expressed thanks to
the lemma \ref{lem_2}, which proves the relation for $S(e_{n-i}(u))$.\\
The proof for the counit is obvious.
\finproof\\
{\it Remark:} This Hopf structure can be extended to the 
double Yangian $DY(sl_n)$. The
image of the generators $x(u)\in Y(sl_n) \subset DY(sl_n)$ 
under the comultiplication, the antipode or the counit are unchanged. 
For the dual generator
$x^*(u)$ of $x(u)$, its image is given by the same formula
where all the generators are replaced by their dual.
\subsection{Examples}
We give two examples to show explicit computations using the
theorems \ref{coproduct} and \ref{th_antipode}.
The comultiplication of $Y(sl_2)$ is given by:

\begin{eqnarray} 
\label{sl2_cop+}
&&\hspace{-20mm}\Delta (e_1(u))=
\sum_{m=0}^{+\infty}\Big(-
e_1(u)\otimes f_1(u+1)
\Big)^m
\Big(
1\otimes e_1(u)+
e_1(u)\otimes \Big(h_1(u)+f_1(u+1)e_1(u)\Big)
\Big)\\
\label{sl2_cop+_mo}
&&=1\otimes e_1(u)
+\sum_{m=0}^{+\infty}(-1)^m
e_1(u)^{m+1}\otimes f_1(u+1)^{m}h_1(u)\\
\label{sl2_cop-}
&&\hspace{-20mm}\Delta (f_1(u))=
\Big(
f_1(u) \otimes 1
+\Big(h_1(u)+f_1(u+1)e_1(u)\Big)\otimes f_1(u)
\Big)
\sum_{m=0}^{+\infty}
\Big(-e_1(u+1)\otimes f_1(u)
\Big)^m
\\
\label{sl2_cop-_mo}
&&=f_1(u)\otimes 1+
\sum_{m=0}^{+\infty}(-1)^m
h_1(u)e_1(u+1)^{m}\otimes f_1(u)^{m+1}\\
\label{sl2_coph}
&&\hspace{-20mm}\Delta (h_1(u))=
\Big(
f_1(u)\otimes e_1(u+1)
+\Big(h_1(u)+f_1(u)e_1(u+1)\Big) \otimes \Big(h_1(u)+f_1(u)e_1(u+1)
\Big)\Big)\nonumber\\
&\times&
\sum_{m=0}^{+\infty}
\Big(-e_1(u+1)\otimes f_1(u)
\Big)^m-\Delta(f_1(u))\Delta(e_1(u+1))
\\
\label{sl2_coph_mo}
&&=\sum_{m=0}^{+\infty}
(-1)^k(k+1)h(u)e_1(u+1)^k\otimes f_1(u+1)^k h(u) 
\end{eqnarray}
The explicit forms (\ref{sl2_cop+_mo}), (\ref{sl2_cop-_mo}) and
(\ref{sl2_coph_mo}) allows us to show that the comultiplication, introduced
in this letter, is the opposite of the comultiplication
used by A.I. Molev \cite{Molev2001}. Remark that the proof of the relation
(\ref{sl2_coph_mo}) from (\ref{sl2_coph}) is not obvious. A
simpler way consists in using
the form (\ref{phih1}) instead of the form (\ref{phih3}) in the proof
of the theorem \ref{coproduct}.
Despite its simpler form, the generalisation to $sl_n$ of the form
given by A.I.Molev does not seem possible.\\
The antipode and the counit of $Y(sl_2)$ are given by:
\begin{eqnarray}
\label{sl2_S+}
S(e_1(u+1))&=&
-e_1(u+1)
\Big(
h_1(u)+f_1(u)e_1(u+1)
\Big)^{-1}\;,\\
\label{sl2_S-}
S(f_1(u+1))&=&
-
\Big(
h_1(u)+f_1(u+1)e_1(u)
\Big)^{-1}
f_1(u+1)\;,\\
\label{sl2_Sh}
S(h_1(u+1))&=&
\Big(
h_1(u)+f_1(u+1)e_1(u)
\Big)^{-1}
-~S(e_1(u+1))~S(f_1(u))\;,
\end{eqnarray}
and 
\begin{eqnarray}
\epsilon(e_1(u))=0~,\quad \epsilon(f_1(u))=0\quad\mbox{and}\quad\epsilon(h_1(u))=1\;.
\end{eqnarray}
For $Y(sl_3)$, the comultiplication 
in the Drinfel'd basis is given by:
\begin{eqnarray}
\label{sl3_cop+}
\Delta (e_1(u))&=&
\sum_{m=0}^{+\infty}(-1)^m\Big(
e_1(u)\otimes f_1(u+1)
-[e_{2}^{0}\;,\;e_1(u)]\otimes [f_{2}^{0}\;,\;f_1(u+1)]
\Big)^m
\\
&&\hspace{-30mm}\times 
\Big(
1\otimes e_1(u)
+ e_1(u)\otimes\Big(h_1(u)+f_1(u+1)e_1(u)\Big)
- [e_{2}^{0}\;,\;e_1(u)]\otimes
\left[f_{2}^{0}\;,h_1(u)+f_1(u+1)e_1(u)\right]
\Big)\nonumber
\end{eqnarray}
\begin{eqnarray}
\label{sl3_cop-}
\Delta(f_1(u))&=&\hspace{-3mm}
\Big(
f_1(u)\otimes 1+
\Big(h_1(u)+f_1(u)e_1(u+1)\Big) \otimes f_1(u)
-[e_{2}^{0}\;,\;h_1(u)+f_1(u)e_1(u+1)]
\otimes[f_{2}^{0}\;,\;f_1(u)]
\Big)\nonumber\\
&&\times
\sum_{m=0}^{+\infty}(-1)^m\Big(
 e_1(u+1)\otimes f_1(u)
-[e_2^{(0)}\;,\;e_1(u+1)] \otimes [f_2^{(0)}\;,\;f_1(u)]
\Big)^m
\\
\label{sl3_coph}
\Delta(h_1(u))&=&
\Big(
f_1(u)\otimes e_1(u)
+
\Big(h_1(u)+f_1(u)e_1(u+1)\Big)
\otimes\Big(h_1(u)+f_1(u)e_1(u+1)\Big)\nonumber\\
&&-
[e_2^{(0)}\;,h_1(u)+f_1(u)e_1(u+1)]
\otimes [f_2^{(0)}\;,h_1(u)+f_1(u)e_1(u+1)]
\Big)\\
&\times&
\sum_{m=0}^{+\infty}(-1)^m\Big(
 e_1(u+1)\otimes f_1(u)
-[e_2^{(0)}\;,\;e_1(u+1)] \otimes [f_2^{(0)}\;,\;f_1(u)]
\Big)^m 
-\Delta(f_1(u))~\Delta(e_1(u+1))\nonumber
\end{eqnarray}
$\Delta(e_2(u))$ (resp. $\Delta(f_2(u))$, $\Delta(h_2(u))$ ) is
obtained by exchanging the subscripts 1 and 2 in 
equation (\ref{sl3_cop+}) (resp. (\ref{sl3_cop-}),
(\ref{sl3_coph})). Of course, the antipode and the counit for
$Y(sl_3)$ can be computed thanks to theorem \ref{th_antipode}, however
we leave it to the reader.

\bigskip
\textbf{Acknowledgements:} I warmly thank D.~Arnaudon, J.~Avan,
V.~Caudrelier, L.~Frappat and E.~Ragoucy for discussions and advice.

\end{document}